\documentclass[reqno,10pt]{amsart}

\usepackage{amsthm,amsmath,amssymb,amsxtra,parskip}
\usepackage[english]{babel}
\usepackage[ansinew]{inputenc} 
\usepackage{graphicx}
\usepackage{subfigure}
\usepackage{geometry}
\usepackage{hyperref}
\usepackage{pst-all}

\newcommand{\C}{\ensuremath{\mathbb{C}}}

\newcommand{\N}{\ensuremath{\mathbb{N}}}

\newcommand{\R}{\ensuremath{\mathbb{R}}}

\renewcommand{\epsilon}{\ensuremath{\varepsilon}}

\newcommand{\Cov}{\text{Cov}}

\newcommand{\Prob}{\mathbb{P}}

\theoremstyle{plain}
\newtheorem{result}{Result}

\theoremstyle{definition}
\newtheorem{definition}{Definition}[]

\theoremstyle{remark}

\begin{document}
\psset{xunit=1cm,yunit=1cm}

\title[Results on random fields with spectral representation]{Some results on random fields admitting a spectral representation with infinitely divisible integrator}
\author{Wolfgang Karcher}
\email{wolfgang.karcher\@@{}uni-ulm.de}
\address{Institute of Stochastics, Ulm University, Helmholtzstr.\ 18, 89069 Ulm, Germany}
\begin{abstract}
We consider random fields $X=\{X(t), t \in \R^d\}$ admitting a spectral representation
$$X(t) = \int_E f_t(x) \Lambda(dx), \quad t \in \R^d,$$
for some set $E$, non-random functions $f_t:E \to \R$, $t \in \R^d$, and an infinitely divisible random measure $\Lambda$ and prove some properties of $X$.
\end{abstract}
\maketitle

\section{Preliminaries}

We start with the definition of an infinitely divisible random measure, cf. \cite{RR89}, pp.~454. Let $E$ be an arbitrary non-empty set and $\mathcal{D}$ be a $\delta$-ring (i.~e. a ring which is closed under countable intersections) of subsets of $E$ such that there exists an increasing sequence $\{E_n\}_{n \in \N} \subset \mathcal{D}$ with $\bigcup_{n \in \N} E_n = E$. Recall that a ring of sets is a non-empty class of sets which is closed under the formation of unions and differences of sets, see e.~g. \cite{Hal74}, p.~19.

Let $\Lambda = \{\Lambda(A): A \in \mathcal{D}\}$ be a real stochastic process defined on some probability space $(\Omega,\mathcal{F}, \Prob)$ such that for each sequence of disjoint sets $\{E_n\}_{n \in \N} \subset \mathcal{D}$, the following properties hold:
\begin{itemize}
 \item $\Lambda$ is \textit{independently scattered}, i.~e. the random variables $\Lambda(E_n), n=1,2,\ldots$, are independent,
 \item $\Lambda$ is \textit{$\sigma$-additive}, i.~e. $\Lambda(\bigcup_{n \in \N} E_n) = \sum_{n \in \N} \Lambda(E_n)$ almost surely if $\bigcup_{n \in \N} E_n \in \mathcal{D}$,
 \item $\Lambda(A)$ is an \textit{infinitely divisible} (\textbf{ID}) random variable for each $A \in \mathcal{D}$, i.~e. $\Lambda(A)$ has the law of the sum of $n$ independent and identically distributed random variables for any natural number $n \in \N$.
\end{itemize}
Then $\Lambda$ is called \textit{infinitely divisible random measure}.

Let $\Psi_{\Lambda(A)}$ be the characteristic function of $\Lambda(A)$. Since $\Lambda(A)$ is \textbf{ID}, its characteristic function is given by the L\'evy-Khintchine representation
\begin{equation}
\Psi_{\Lambda(A)}(t) = \exp\left\{it\nu_0(A) - \frac{1}{2} t^2\nu_1(A) + \int_{\R}\left(e^{itx} - 1- it\tau(x)\right)F_A(dx)\right\}, \label{eq:cf_Lambda}
\end{equation}
where $\nu_0: \mathcal{D} \to \R$ is a signed measure, $\nu_1: \mathcal{D} \to [0,\infty)$ is a measure, $F_A:\R \to [0,\infty)$ is a L\'evy measure, i.~e.
$$\int_\R \min(1,z^2)F_A(dz) < \infty$$
and
$$ \tau(z) = \begin{cases}
			  z,& \vert z \vert \leq 1, \\
			  \frac{z}{\vert z \vert}, & \vert z \vert > 1.
			 \end{cases}$$
Define the measure $\lambda$ by
$$\lambda(A) := \vert \nu_0 \vert(A) + \nu_1(A) + \int_\R \min(1,z^2) F_A(dz), \quad A \in \mathcal{D}.$$
We call $\lambda$ \textit{control measure} of the \textbf{ID} random measure $\Lambda$.

Let $\sigma(\mathcal{D})$ be the $\sigma$-algebra generated by $\mathcal{D}$ and $I_A: E \to \{0,1\}$ the indicator function of a set $A \subset E$ with
\begin{equation*}
 I_A(x) := \begin{cases}
			1, & x \in A,\\
			0, & x \notin A.
		   \end{cases} \label{eq:ind}
\end{equation*}

For disjoint sets $A_j \in \mathcal{D}$, real numbers $x_j$, $j=1,\ldots,n$, $n \in \N$, and simple functions of the form $f = \sum_{j=1}^n x_j I_{A_j}$, we define for every $A \in \sigma(\mathcal{D})$
$$\int_A f d\Lambda := \sum_{j=1}^n x_j \Lambda(A \cap A_j).$$

Let $f_t:E \to \R$, $t \in \R^d$, $d \in \N$, be a $\sigma(\mathcal{D})$-measurable function which is \textit{$\Lambda$-integrable}, that is there exists a sequence of simple functions $\{\tilde{f}_t^{(n)}\}_{n \in \N}$ such that
\begin{enumerate}
 \item $\tilde{f}_t^{(n)} \to f_t \quad \lambda-\text{a.e.},$
 \item for every set $A \in \sigma(\mathcal{D})$, the sequence $\{\int\limits_{A} \tilde{f}_t^{(n)}(x) \Lambda(dx)\}_{n \in \N}$ converges in probability as $n \to \infty$.
\end{enumerate}

A family $X = \{X(t), t \in \R^d\}$ of real-valued random variables $X(t)$ defined on a probability space $(\Omega,\mathcal{F},\Prob)$ is called \textit{random field}. For each $t \in \R^d$, we define
\begin{equation}
 \int\limits_{E} f_t(x) \Lambda(dx) := \underset{n \to \infty}{\text{plim}} \int\limits_{E} \tilde{f}_t^{(n)}(x) \Lambda(dx), \label{eq:def_int}
\end{equation}
where $plim$ means convergence in probability, and consider random fields of the form
\begin{equation}
 X(t) = \int\limits_{E} f_t(x) \Lambda(dx), \quad t \in \R^d. \label{eq:spectralRepresentation}
\end{equation}
In \cite{UW67}, it is shown that (\ref{eq:def_int}) does not depend on the approximation sequence $\{\tilde{f_t}^{(n)}\}_{n \in \N}$ and thus is well-defined.

Notice that by Lemma~2.3. in \cite{RR89}, we have
\begin{eqnarray*}
 F_A(B) = F(A \times B) \quad \text{and} \quad F(dx,ds) = \rho(x,ds)\lambda(dx),
\end{eqnarray*}
where $F$ is a $\sigma$-finite measure on $\sigma(\mathcal{D}) \times \mathcal{B}(\R)$. Furthermore, $\rho: E \times \mathcal{B}(\R) \to [0,\infty]$ is a function such that $\rho(x,\cdot)$ is a L\'evy measure on $\mathcal{B}(\R)$ for every $x \in E$ and $\rho(\cdot,B)$ is a Borel measurable function for every $B \in \mathcal{B}(\R)$. Moreover, $\nu_0$ and $\nu_1$ are absolutely continuous with respect to $\lambda$. We set $a := d\nu_0 / d\lambda$ and $\sigma^2 := d\nu_1 / d\lambda$.

Let us introduce a certain type of dependence structure, namely (positive or negative) association, which we will consider in the following section. Let $\vert I \vert$ denote the cardinality of a finite set $I \subset T$ and $X_I:= \{X(t), t \in I\}$.

\begin{definition}\label{def:association}
 Let $\mathcal{M}(n)$ be the class of real-valued bounded coordinate-wise nondecreasing Borel functions on $\R^n$, $n \in \N$, and $T$ be an index set.

\begin{enumerate}
 \item[(a)] A family $\{X(t), t \in T\}$ is called associated (\textbf{A}) if for every finite set $I \subset T$ and any functions $f,g \in \mathcal{M}(\vert I\vert)$, one has
 $$ Cov(f(X_I),g(X_I)) \geq 0.$$
 \item[(b)] A family $\{X(t), t \in T\}$ is called positively associated (\textbf{PA}) if for any disjoint finite sets $I,J \subset T$ and all functions $f \in \mathcal{M}(\vert I\vert)$, $g \in \mathcal{M}(\vert J\vert)$, one has
 $$ Cov(f(X_I),g(X_J)) \geq 0.$$
 \item[(c)] A family $\{X(t), t \in T\}$ is called negatively associated (\textbf{NA}) if for any disjoint finite sets $I,J \subset T$ and all functions $f \in \mathcal{M}(\vert I\vert)$, $g \in \mathcal{M}(\vert J\vert)$, one has
 $$ Cov(f(X_I),g(X_J)) \leq 0.$$
\end{enumerate}
\end{definition}

In the above definition, any permutation of coordinates of the random vector $(X(t_1),\ldots,X(t_n))^\mathsf{T}$ is used for $X_K$, $K=\{t_1,\ldots,t_n\} \subset T$.

Finally, we provide the definition of stochastic continuity.

\begin{definition}
A random field $X=\{X(t), t \in \R^d\}$ is called \textit{stochastically continuous} at $t \in \R^d$ if $\underset{s \to t}{\text{plim}}X(s) = X(t)$.
\end{definition}

\section{Results}

\begin{result}
 Let $X=\{X(t), t \in \R^d\}$ be a random field of the form (\ref{eq:spectralRepresentation}). Then $X$ is \textbf{ID}, that is the law of the random vector $(X(t_1),\ldots,X(t_n))^\mathsf{T}$, $n \in \N$, is an \textbf{ID} probability measure on $\R^n$ for all $t_1,\ldots,t_n \in \R^d$.
\end{result}

\begin{proof}
Let $\varphi_{(t_1,\ldots,t_n)}$ be the characteristic function of $(X(t_1),\ldots,X(t_n))^\mathsf{T}$. It is enough to show that $\varphi_{(t_1,\ldots,t_n)}^\gamma$ is a characteristic function for all $\gamma > 0$.

 Due to the linearity of the spectral representation (\ref{eq:spectralRepresentation}) and the fact that any linear combination of $\Lambda$-integrable functions is $\Lambda$-integrable (cf. \cite{JW94}, p.~81), we have
  \begin{eqnarray*}
  \varphi_{(t_1,\ldots,t_n)}^\gamma(x) &=& \varphi_{\sum_{j=1}^n x_j X(t_j)}(1) \\
  &=& \exp\left\{\int_E \left[\sum_{j=1}^n x_j f_{t_j}(y) a(y) - \frac{1}{2} t^2 \left(\sum_{j=1}^n x_j f_{t_j}(y)\right)^2 \sigma^2(y)\right.\right.\\
 &&\left.\left. + \int_\R \left(e^{it \sum_{j=1}^n x_j f_{t_j}(y) s} - 1 - it \sum_{j=1}^n x_j f_{t_j}(y) \tau(s) \right) \rho(y,ds)\right] \gamma \lambda(dy)\right\},
  \end{eqnarray*}
 where the last equality follows from Proposition~2.6. in \cite{RR89}. Define $\nu_0^*:\mathcal{D} \to \R$, $\nu_1^*: \mathcal{D} \to [0,\infty)$, $F_A^*:\mathcal{B}(\R) \to [0,\infty)$ by $\nu_0^*(ds) := a(s) \gamma \lambda(ds) = \gamma \nu_0(ds)$, $\nu_1^*(ds):= \sigma^2(s) \gamma \lambda(ds) = \gamma \nu_1(ds)$ and
 $$F_A^*(B) := \int_E \int_\R I_{A \times B}(s,x) \rho(s,dx) \gamma \lambda(ds) = \gamma \int_{A \times B} F(ds,dx) = \gamma F(A \times B) = \gamma F_A(B)$$
 for all $A \in \mathcal{D}$ and $B \in \mathcal{B}(\R)$, cf. Lemma~2.3. in \cite{RR89}. Since $\nu_0^*$ is a signed measure, $\nu_1^*$ is a measure, $F_A^*$ is a L\'evy measure on $\R$ for all $A \in \mathcal{D}$ and $F.(B)$ is a measure for all $B \in \mathcal{B}(\R)$ whenever $0 \notin \bar{B}$, there exists an \textbf{ID} random measure $\Lambda^*$ with characteristic function (\ref{eq:cf_Lambda}) (and control measure $\lambda^*=\gamma \lambda$), where $\nu_0$, $\nu_1$ and $F_A$ in (\ref{eq:cf_Lambda}) is replaced by $\nu_0^*$, $\nu_1^*$ and $F_A^*$, respectively, see Proposition~2.1.~(b) in \cite{RR89}. Therefore, $\varphi_{(t_1,\ldots,t_n)}^\gamma$ is the characteristic function of $(Y(t_1),\ldots,Y(t_n))^\mathsf{T}$ with
$$Y(t) := \int_E f_t(x) \Lambda^*(dx).$$
\end{proof}

The following result provides conditions for $\Lambda$-integrability, cf. Lemma 1 in \cite{HPVJ08}.

\begin{result}\label{result:integrability}
 Let $f:E \to \R$ be a $\sigma(\mathcal{D})$-measurable function. If
 \begin{itemize}
  \item[(i)] $\int_E \vert f(x) a(x) \vert \lambda(dx) < \infty$,
  \item[(ii)] $\int_E f^2(x) \sigma^2(x) \lambda(dx) < \infty$,
  \item[(iii)] $\int_E \int_\R \vert f(x) s \vert \rho(x,ds) \lambda(dx) < \infty$,
 \end{itemize}
 then $f$ is $\Lambda$-integrable and the characteristic function of $\int_E f(x) \Lambda(dx)$ is given by
 \begin{eqnarray*}
 &&\hspace*{-0.5cm}\Psi_{\int_E f(x) \Lambda(dx)}(t) \\
  &&\hspace*{-0.5cm}= \exp\left\{ it \int_E f(x) \nu_0(dx) - \frac{1}{2} t^2 \int_E f^2(x) \nu_1(dx) + \int_E \int_\R \left(e^{itf(x) s}-1-itf(x) \tau(s)\right) F(dx,ds)\right\}.
 \end{eqnarray*}
\end{result}

\begin{proof}
By Theorem~2.7 in \cite{RR89}, it suffices to show
\begin{itemize}
 \item[(a)] $\int_E \vert U(f(x),x)\vert \lambda(dx) < \infty$,
 \item[(b)] $\int_E \vert V_0(f(x),x)\vert \lambda(dx) < \infty$,
\end{itemize}
where
\begin{eqnarray*}
U(u,x) &=& u a(x) + \int_\R \left(\tau(su) - u\tau(s)\right)\rho(x,ds) \\
V_0(u,x) &=& \int_\R \min\{1,\vert su\vert^2\} \rho(x,ds).
\end{eqnarray*}
We follow the proof of Lemma~1 in \cite{HPVJ08}. It holds $\vert \tau(su) \vert \leq \vert su \vert$. This implies
\begin{equation*}
 \vert U(f(x),x)\vert \leq \vert f(x) a(x) \vert + \int_\R \left(\vert f(x)s \vert + \vert f(x)s\vert \right)\rho(x,ds) = \vert f(x) a(x) \vert + 2 \int_\R \vert f(x)s\vert\rho(x,ds)
\end{equation*}
such that condition (a) is satisfied by (i) and (iii). Since $\min\{1,(s f(x))^2\} \leq \vert s f(x) \vert$, condition (b) is satisfied by (iii).

We now derive the formula for the characteristic function of $\int_E f(x) \Lambda(dx)$. By Proposition 2.6. in \cite{RR89}, it is given by
\begin{eqnarray*}
 &&\hspace*{-0.5cm}\Psi_{\int_E f(x) \Lambda(dx)}(t) \\
 &&\hspace*{-0.5cm}= \exp\left\{ it \int_E \left[f(x) a(x) - \frac{1}{2} t^2 f^2(x) \sigma^2(x) + \int_\R \left(e^{itf(x) s}-1-itf(x) \tau(s)\right) \rho(x,ds)\right] \lambda(dx)\right\}.
\end{eqnarray*}
We have $\int_E \vert f(x) a(x) \vert \lambda(dx) < \infty$ and $\int_E f^2(x) \sigma^2(x) \lambda(dx) < \infty$ by (i) and (ii). It remains to show that
$$\int_E \left\vert \int_\R \left(e^{itf(x) s}-1-itf(x) \tau(s)\right) \rho(x,ds) \right\vert \lambda(dx) < \infty.$$

Let $y \neq 0$. By using the mean value theorem, we get
\begin{eqnarray*}
 \left\vert \frac{\sin(y) - \sin(0)}{y-0} \right\vert &=& \vert \sin(\xi_1)\vert \leq 1, \\
 \left\vert \frac{\cos(y) - \cos(0)}{y-0} \right\vert &=& \vert \cos(\xi_2)\vert \leq 1,
\end{eqnarray*}
where $\xi_1,\xi_2 \in [0,y]$ if $y>0$ and $\xi_1,\xi_2 \in [y,0]$ if $y<0$. Therefore, we have for each $y \in \R$
\begin{eqnarray*}
 \vert e^{iy} - 1\vert &\leq& \vert e^{iy} - e^{i0} \vert = \vert \cos(y) + i \sin(y) - \cos(0) - i \sin(0) \vert \\
 &=& \sqrt{(\cos(y) - \cos(0))^2+(\sin(y)-\sin(0))^2} \leq \sqrt{y^2 + y^2} = \sqrt{2} \vert y \vert.
\end{eqnarray*}
This implies
\begin{eqnarray*}
&&\int_E \left\vert \int_\R \left(e^{itf(x) s}-1-itf(x) \tau(s)\right) \rho(x,ds) \right\vert \lambda(dx) \\
&\leq& \int_E\int_\R\left(\left\vert e^{itf(x) s}-1\right\vert + \vert tf(x) \tau(s)\vert\right) \rho(x,ds)\lambda(dx) \\
&\leq& \int_E\int_\R\left( \sqrt{2} \vert t f(x) s \vert + \vert t f(x) s\vert \right) \rho(x,ds)\lambda(dx) \\
&\leq& t (\sqrt{2}+1) \int_E\int_\R \vert f(x) s \vert \rho(x,ds)\lambda(dx) < \infty,
\end{eqnarray*}
where the last inequality follows from (iii).
\end{proof}

We now provide a sufficient condition for the independence of two families of random variables taken from the random field (\ref{eq:spectralRepresentation}). We denote the support of a function $f$ by $supp(f)$.

\begin{result}\label{lemma:IDindependence}
Let $X$ be a random field of the form (\ref{eq:spectralRepresentation}). Let $K,L \subset T=\{t_1,\ldots,t_k\}$, $T \subset \R^d$, $k\in\N$, with $K\cup L = T$, $K,L \neq \emptyset$ and $K\cap L = \emptyset$. If
\begin{equation}
\left(\bigcup_{t_i \in K} supp(f_{t_i})\right)\bigcap\left(\bigcup_{t_j \in L} supp(f_{t_j})\right) = \emptyset, \label{eq:support}
\end{equation}
then the families of random variables $\{X(t_i), t_i \in K\}$ and $\{X(t_j), t_j \in L\}$ are independent.
\end{result}

\begin{proof}
 Let $\varphi_K$, $\varphi_L$ and $\varphi_{T}$ be the characteristic functions of a fixed permutation of the random vectors constructed from the families of random variables $\{X(t_i), t_i \in K\}$, $\{X(t_j), t_j \in L\}$ and $\{X(t), t \in T\}$. Furthermore, let $x_K \in \R^{\vert K \vert}$, $x_L \in \R^{\vert L \vert}$ and $x_T \in \R^{\vert T \vert}$ and define $K:\R\times E \to \C$ by
 $$K(t,s) := it a(s) - \frac{1}{2} t^2 \sigma^2(s) + \int_\R\left(e^{itx} - 1 - it\tau(x)\right) \rho(s,dx).$$ 
We have
 \begin{eqnarray*}
 \varphi_{K} (x_K) &=& \varphi_{\sum\limits_{t_i \in K} x_{t_i} X(t_i)}(1) = \exp\left\{\int_E K\left(\sum_{t_i \in K} x_{t_i} f_{t_i}(s),s\right) \lambda(ds)\right\}, \\
 \varphi_{L} (x_L) &=& \exp\left\{\int_E K\left(\sum_{t_j \in L} x_{t_j} f_{t_j}(s),s\right) \lambda(ds)\right\}, \\
 \varphi_{T} (x_T) &=& \exp\left\{\int_E K\left(\sum_{t \in T} x_{t} f_{t}(s),s\right) \lambda(ds)\right\},
 \end{eqnarray*}
 see Proposition~2.6. in \cite{RR89}. By using condition (\ref{eq:support}), it is not difficult to check that for each $s \in \cup_{t_i \in K} supp(f_{t_i})$ and each $s \in \cup_{t_j \in L} supp(f_{t_j})$
 \begin{equation}
 K\left(\sum_{t_i \in K} x_{t_i} f_{t_i}(s),s\right) + K\left(\sum_{t_j \in L} x_{t_j} f_{t_j}(s),s\right) = K\left(\sum_{t \in T} x_{t} f_{t}(s),s\right). \label{eq:help28}
 \end{equation}
 Thus, (\ref{eq:help28}) holds for all $s \in E$. This implies
 $$\varphi_K(x_K) \varphi_L(x_L) = \varphi_T(x_T) $$
 such that $\{X(t_i), t_i \in K\}$ and $\{X(t_j), t_j \in L\}$ are independent (cf.~Theorem~4 in \cite{Shi96}, p.~286, and its proof).
\end{proof}

The following result provides a sufficient condition for a random field of the form (\ref{eq:spectralRepresentation}) to be associated.
\begin{result}\label{result:association}
Suppose that for all $t \in \R^d$, either $f_t(x) \geq 0$  for all $x \in E$ or $f_t(x) \leq 0$ for all $x \in E$. Then (\ref{eq:spectralRepresentation}) is an associated random field.
\end{result}

\begin{proof}
Let $f$ be $\Lambda$-integrable and non-negative. In the proof of Theorem~2.7. in \cite{RR89}, the corresponding approximating sequence $\{\tilde{f}^{(n)}\}_{n \in \N}$ for $f$ is selected in the following way.

Let $A_n = \{x \in E: \vert f(x) \vert \leq n \} \cap E_n$. Choose a sequence of simple $\mathcal{D}$-measurable functions $\{\tilde{f}^{(n)}\}_{n \in \N}$ such that
\begin{equation}
f^{(n)}(x) = 0 \text{ if }x \notin A_n, \quad \vert f^{(n)}(x) - f(x) \vert \leq \frac{1}{n} \text{ if } x \in A_n, \quad \ \vert f^{(n)}(x) \vert \leq \vert f(x) \vert\,\, \forall x \in E. \label{eq:approx_props}
\end{equation}
We now define a simple function $f_*^{(n)}$ by
$$f_*^{(n)}(x):=\begin{cases}
	f^{(n)}(x), & f_n(x) \geq 0, \\
	0, & f_n(x) < 0,
  \end{cases} \quad x \in E.$$
Since $f$ is non-negative, it is easy to see that the sequence $\{f_*^{(n)}\}_{n \in \N}$ fulfills the same properties (\ref{eq:approx_props}) as $\{\tilde{f}^{(n)}\}_{n \in \N}$. So $\{f_*^{(n)}\}_{n \in \N}$ is an approximating sequence for $f$ which is additionally non-negative. As $f_*^{(n)}$ is simple for all $n \in \N$, we can write
$$f_*^{(n)} = \sum_{j=1}^{m(n)} x_j I_{B_j}$$
for some $m(n) \in \N$, $x_j \geq 0$ and disjoint $B_j \subset A_n$, $j=1,\ldots,m(n)$.

Assume now that for all $t \in \R^d$, $f_t(x) \geq 0$ for all $x \in E$. Let $t_1,\ldots,t_r \in \R^d$, $r \in \N$, and $\{f_{t_i}^{(n)}\}_{n \in \N}$ be the approximating sequences of the kernel functions $f_{t_i}$ in the spectral representation
$$X(t_i) = \int_E f_{t_i} \Lambda(dx).$$
Consider
$$X^{(n)}(t_i) = \int_E f_{t_i}^{(n)}(x) \Lambda(dx) = \int_E \sum_{j=1}^{m(n,i)} x_j^{(i)} I_{B_j^{(i)}}(x)\Lambda(dx) = \sum_{j=1}^{m(n,i)} x_j^{(i)} \Lambda(B_j^{(i)})$$
and, as we have just seen, we can assume without loss of generality that $x_j^{(i)} \geq 0$. By further decomposing $B_j^{(i)}$ if necessary, we can find a set of disjoint sets $\{B_j,j=1,\ldots,m(n)\}$ for some $m(n) \in \N$ which does not depend on $i$ such that
$$X^{(n)}(t_i) = \sum_{j=1}^{m(n)} x_j^{(i)} \Lambda(B_j), \quad \forall i=1,\ldots,r.$$
The sets $B_j$, $j=1,\ldots,m(n)$, can be obtained by intersecting and subtracting the sets $B_j^{(i)}$ appropriately. This implies that $B_j \in \mathcal{D}$, $j=1,\ldots,m(n)$, by using the properties of rings of sets.

We now show that the random vector $(X^{(n)}(t_1),\ldots,X^{(n)}(t_r))^\mathsf{T}$ is \textbf{A} for all $n \in \N$. Consider the set $I=\{t_1,\ldots,t_r\}$ and let $f \in \mathcal{M}(r)$ and $g \in \mathcal{M}(r)$. We write $X^{(n)}(I)$ for the random vector consisting of an arbitrary permutation of the corresponding elements $X^{(n)}(s)$, $s \in I$. Consider the functions $k,l:\R^{m(n)} \to \R$ defined by
\begin{eqnarray*}
k(y_1,\ldots,y_{m(n)}) &:=& f\left(\sum_{j=1}^{m(n)} x_j^{(1)} y_j,\ldots,\sum_{j=1}^{m(n)} x_j^{(r)} y_j\right), \quad (y_1,\ldots,y_{m(n)})^\mathsf{T} \in \R^{m(n)}, \\
l(y_1,\ldots,y_{m(n)}) &:=& g\left(\sum_{j=1}^{m(n)} x_j^{(1)} y_j,\ldots,\sum_{j=1}^{m(n)} x_j^{(r)} y_j\right), \quad (y_1,\ldots,y_{m(n)})^\mathsf{T} \in \R^{m(n)}.
\end{eqnarray*}
Since the coefficients $x_j^{(i)}$ are non-negative for all $i=1,\ldots,r$ and $f \in \mathcal{M}(r)$, $g \in \mathcal{M}(r)$, we conclude that $k,l \in \mathcal{M}(m(n))$.

By definition, $\Lambda(B_1),\ldots,\Lambda(B_{m(n)})$ are independent and therefore \textbf{A}, cf. Theorem 1.8~(c) in \cite{BS07}, p.~6. This implies
$$\Cov\left(f(X^{(n)}(I)),g(X^{(n)}(I))\right) = \Cov\left(k(\Lambda(B_1),\ldots,\Lambda(B_{m(n)})),l(\Lambda(B_1),\ldots,\Lambda(B_{m(n)}))\right) \geq 0 $$
such that $(X^{(n)}(t_1),\ldots,X^{(n)}(t_r))^\mathsf{T}$ is associated.

For a vector $x=(x_1,\ldots,x_r)^\mathsf{T} \in \R^r$, set $\Vert x \Vert_1 := \sum_{i=1}^r \vert x_i \vert$. Since $X^{(n)}(t_i)$ converges in probability to $X(t_i)$ as $n \to \infty$ for $i=1,\ldots,r$, we conclude that $\Vert (X^{(n)}(t_1),\ldots,X^{(n)}(t_r))^\mathsf{T} \Vert_1$ converges to $\Vert (X(t_1),\ldots,X(t_r))^\mathsf{T} \Vert_1$ in probability due to Markov's inequality. This implies that $(X^{(n)}(t_1),\ldots,X^{(n)}(t_r))^\mathsf{T}$ converges in distribution to $(X(t_1),\ldots,X(t_r))^\mathsf{T}$, see \cite{Bil68}, p.~18. Therefore, $(X(t_1),\ldots,X(t_r))^\mathsf{T}$ is \textbf{A}, see \cite{BS07}, p.~7. Thus $X$ is an associated random field.

If for all $t \in \R^d$, the functions $f_t$ are non-positive for all $x \in E$, the proof is completely analogous by considering the fact that in Definition~\ref{def:association} (a) one can use coordinate-wise non-increasing functions instead of coordinate-wise non-decreasing functions, see Remark~1.4. in \cite{BS07}, p.~4. In this case, the approximating sequences $\{f_{t_i}^{(n)}\}_{n \in \N}$ of the kernel functions $f_{t_i}$, $i=1,\ldots,r$, are chosen in such a way that $f_{t_i}^{(n)}$ is non-positive for each $n \in \N$.
\end{proof}

The following result provides sufficient conditions such that a random field of the form (\ref{eq:spectralRepresentation}) is stochastically continuous.

\begin{result}\label{lemma:stochCont}
Assume that the following conditions for a random field $X$ with spectral representation (\ref{eq:spectralRepresentation}) hold. For each $t \in \R^d$,
\begin{itemize}
\item[(a)] $f_s \to f_t$ $\lambda$-almost everywhere as $s \to t$,
\item[(b)] there exists some $\varepsilon > 0$ and a $\Lambda$-integrable function $g$ such that $\vert f_s - f_t \vert \leq g$ $\lambda$-almost everywhere and for all $s \in \R^d$ such that $\Vert s-t \Vert_2 \leq \varepsilon$,
\end{itemize}
where $\lambda$ is the control measure of the \textbf{ID} random measure $\Lambda$. Then $X$ is stochastically continuous.
\end{result}

\begin{proof}
 In \cite{RR89}, the discussion before Theorem~3.3. and Theorem~3.3. itself imply that there is a function $\Phi_0:\R \times E \to [0,\infty)$ such that for a sequence of $\Lambda$-measurable functions $\{f_n\}_{n \in \N}$ we have the following implication:
 \begin{equation}
 \int_E \Phi_0(\vert f_n(x) \vert, x)\lambda(dx) \to 0, \quad n \to \infty \quad \Rightarrow\quad \underset{n \to \infty}{\text{plim}}\int_E f_n(x) \Lambda(dx) = 0.\label{eq:help26}
 \end{equation}
 Furthermore for every $x \in E$, $\Phi_0(\cdot,x)$ is a continuous non-decreasing function on $[0,\infty)$ with $\Phi_0(0,x)=0$, see Lemma~3.1. in \cite{RR89}. Since $g$ in assumption (b) is $\Lambda$-integrable, we have
 $$\int_E \Phi_0(\vert g(x) \vert, x)\lambda(dx) < \infty,$$
 cf. the definition of the Musielak-Orlicz space $L_{\Phi_0}(E,\lambda)$ on p.~466 and again Theorem~3.3 in \cite{RR89}.
 
 Let $t \in \R^d$. As $\Phi_0(\cdot,x)$ is non-decreasing, we get
 $$\int_E \Phi_0(\vert f_s(x)-f_t(x) \vert, x)\lambda(dx) \leq \int_E \Phi_0(\vert g(x) \vert, x)\lambda(dx)$$
 for all $s \in \R^d$ such that $\Vert s-t \Vert_2 \leq \varepsilon$ by assumption (b). Furthermore
 $$ \Phi_0(\vert f_s(x)-f_t(x) \vert, x) \to \Phi_0(0,x) = 0 \quad \lambda-a.e.$$
 as $s \to t$ due to the continuity of $\Phi_0(\cdot,x)$ and assumption (a). Therefore, we can apply the dominated convergence theorem and get
 $$\int_E \Phi_0(\vert f_s(x)-f_t(x) \vert, x)\lambda(dx) \to 0, \quad s \to t.$$
 Then (\ref{eq:help26}) implies
 $$\underset{s \to t}{\text{plim}} \int_E (f_s(x)-f_t(x))\Lambda(dx) = \underset{s \to t}{\text{plim}} \left(\int_E f_s(x)\Lambda(dx) - \int_E f_t(x)\Lambda(dx)\right) = 0,$$
 that is
 $$\underset{s \to t}{\text{plim}} X(s) = \underset{s \to t}{\text{plim}}\int_E f_s(x)\Lambda(dx) = \int_E f_t(x)\Lambda(dx) = X(t).$$
\end{proof}

If $\Lambda=M$ is an $\alpha$-stable random measure (cf. \cite{ST94}, pp. 118), then
$$X(t) = \int_E f_t(x) M(dx), \quad t \in \R^d,$$
is an $\alpha$-stable random field since the random vector $(X(t_1),\ldots,X(t_n))^\mathsf{T}$ is multivariate $\alpha$-stable distributed  for all $t_1,\ldots,t_n \in \R^d$, $n \in \N$, cf. Proposition 3.4.3 in \cite{ST94}, p.~125.

Recall that the characteristic function of a stable random vector $\text{\boldmath{$X$}}=(X_1,\ldots,X_n)^\mathsf{T}$, $n \in \N$, is given by
\begin{eqnarray*}
 \varphi_{\text{\boldmath{$X$}}}(\text{\boldmath{$\theta$}})=E\left(e^{i \cdot \text{\boldmath{$\theta$}}^\mathsf{T} \text{\boldmath{$X$}}}\right) =\begin{cases}
               e^{-\int_{S_n}|\text{\boldmath{$\theta$}}^\mathsf{T} \text{\boldmath{$s$}}|^{\alpha} (1-i(\text{sign}\,
\text{\boldmath{$\theta$}}^\mathsf{T} \text{\boldmath{$s$}})\tan \frac{\pi \alpha}{2}\Gamma(d\text{\boldmath{$s$}}) + i
\text{\boldmath{$\theta$}}^\mathsf{T} \text{\boldmath{$\mu$}} } \quad \text{if } \alpha \ne 1,\\
               e^{-\int_{S_n}|\text{\boldmath{$\theta$}}^\mathsf{T} \text{\boldmath{$s$}}| (1+i\frac{2}{\pi}(\text{sign}\,
\text{\boldmath{$\theta$}}^\mathsf{T} \text{\boldmath{$s$}})\ln |\text{\boldmath{$\theta$}}^\mathsf{T}
\text{\boldmath{$s$}}|\Gamma(d\text{\boldmath{$s$}}) + i \text{\boldmath{$\theta$}}^\mathsf{T} \text{\boldmath{$\mu$}} } \quad
\text{\hspace*{-0.23cm} if } \alpha = 1, \\
             \end{cases}\quad \forall \text{\boldmath{$\theta$}} \in \R^n,
\end{eqnarray*}
where $\Gamma$ is a finite measure on the unit sphere $S_n$ of $\R^n$ and $\text{\boldmath{$\mu$}} \in \R^n$.

An $\alpha$-stable random vector $\text{\boldmath{$X$}}=(X_1,\ldots,X_n)^\mathsf{T}$, $n \in \N$, is \textbf{A} if and only if $\Gamma(S_-) = 0$,
where $S_- := \{ (s_1,\ldots,s_n) \in S_n: s_i s_j < 0 \text{ for some } i,j\}$. It is \textbf{NA} if and only if $\Gamma(S_+) = 0$, where $S_+ := \{ (s_1,\ldots,s_n) \in
S_n: s_i s_j > 0 \text{ for some } i,j\}$, see Theorem 4.6.1, p.~204, and Theorem 4.6.3, p.~208, in \cite{ST94}. The following result yields the same sufficient condition for association as Result \ref{result:association}, but provides a more straightforward proof.

\begin{result} \label{result:posNeg}
 Suppose that for all $t \in \R^d$, either $f_t(x) \geq 0$  for all $x \in E$ or $f_t(x) \leq 0$ for all $x \in E$. Then
 $$  X(t) = \int_E f_t(x) M(dx), \quad t \in \R^d,$$
 is an associated $\alpha$-stable random field.
\end{result}

\begin{proof}
 Consider the finite random vector $(X(t_1),\ldots,X(t_n))^\mathsf{T}$, $n \in \N$, and define the set $E_+$ and the function
$g=(g_1,\ldots,g_n):E_+ \to \R^n$ by
 \begin{eqnarray*}
  E_+ &:=& \{x \in E: \sum_{k=1}^n f_{t_k}(x)^2 > 0\}, \\
  g_j(x) &:=& \frac{f_{t_j}(x)}{(\sum_{k=1}^n f_{t_k}(x)^2)^{1/2}}, \quad j=1,\ldots,n.
 \end{eqnarray*}
 Then, for any Borel set $A$ in $S_n$, we have with $-A:=\{-a: a \in A\}$
 $$\Gamma(A) = \int_{g^{-1}(A)} \frac{1+\beta(x)}{2} m_1(dx) + \int_{g^{-1}(-A)} \frac{1-\beta(x)}{2} m_1(dx),$$
 where
 \begin{eqnarray*}
  m_1(dx) &=& \left(\sum_{k=1}^n f_{t_k}(x)^2\right)^{\alpha/2} m(dx), \\
  g^{-1}(A) &=& \{x \in E_+: (g_1(x),\ldots,g_n(x)) \in A\},
 \end{eqnarray*}
 see \cite{ST94}, pp.~115. Since $g_i(x) g_j(x) \geq 0$ for all $x \in E_+$ and $i,j \in \{1,\ldots,n\}$, we get
$g^{-1}(S_-) = \emptyset$ and $g^{-1}(-S_-) = \emptyset$ and therefore $\Gamma(S_-) = 0$. Thus, $X$ is \textbf{A}.
\end{proof}

\begin{result}\label{result:as_null}
Let $M$ be an $\alpha$-stable random measure with control measure $m$ and let $f$ be $M$-integrable. If $\int_E \vert f(x)\vert^\alpha m(dx) = 0$, then $\int_E f(x) M(dx) = 0$ almost surely.
\end{result}

\begin{proof}
Let $g \in F$, where $F$ is the set of all $M$-integrable functions. Since $\alpha$-stable integrals are linear (see~\cite{ST94}, p.~125), we have
$$\int_E 0 M(dx) = \int_E 0 g(x) M(dx) = 0 \int_E g(x) M(dx) = 0 \quad a.s.$$
Notice also that for the null function $h:E \to \R$ with $h(x):=0$ for all $x \in E$, we have $h \in F$ since $F$ is a linear space (see \cite{ST94}, p.~122). The assumption $\int_E \vert f(x)\vert^\alpha m(dx) = 0$ implies that $f=0$ $m$-almost everywhere. We can therefore use $\{f^{(n)}\}_{n \in \N}$ with $f^{(n)} = h$ as an approximating sequence for $f$ which has the properties (3.4.7) and (3.4.8) in \cite{ST94}, p.~122. We have
$$\underset{n \to \infty}{\text{plim}}\int_E f^{(n)}(x) M(dx) = \underset{n \to \infty}{\text{plim}} 0 = 0,$$
but also
$$\underset{n \to \infty}{\text{plim}}\int_E f^{(n)}(x) M(dx) = \int_E f(x) M(dx),$$
see \cite{ST94}, p.~124. Since convergence in probability implies convergence in distribution and the corresponding limit distribution is unique (see \cite{Bil68}, p.~11 and p.~18), the result is proven.
\end{proof}

\end{document}